# A Direct Method To Generate Pythagorean Triples And Its Generalization To Pythagorean Quadruples And n-tuples

*Tanay Roy$ and Farjana Jaishmin Sonia£*
Department of Physics
Jadavpur University
Kolkata 700032, India
$roytanayslg@gmail.com, £farjana.jas@gmail.com

**Abstract:** The method of generating Pythagorean triples is known for about 2000 years. Though the classical formulas produce all primitive triples they do not generate all possible triples, especially non-primitive triples. This paper presents a direct method to generate all possible triples both primitive and non-primitive for any given number. Then the technique is generalized to produce Pythagorean Quadruples and n-tuples. Our method utilizes the fact that the difference between lengths of the hypotenuse and one leg of a Pythagorean triangle can have only certain distinct values depending on the length of the other which remains true for higher tuples also.

## I. Introduction

A Pythagorean triple is an ordered triple of positive integers $(a, b, c)$ such that
$$a^2 + b^2 = c^2 \qquad (1)$$
An effective way to generate Pythagorean triples is based on Euclid's formula found in his book *Elements*. This formula states that for any two positive integers $m$ and $n$ with $m > n$ $a = m^2 - n^2$, $b = 2mn$, $c = m^2 + n^2$ form a Pythagorean triple. Though this classical formula generates all primitive triples and infinitely many of them, it is incapable of generating all the triples including non-primitive cases. For example the triple (15,36,39) cannot be generated from the formula rather introducing a multiplier to the triple (5,12,13) does so. Besides we observe that while Euclid's formula produces the triple (5,12,13), it doesn't produce (12,5,13); a transposition is needed. Also by a result of Berggren (1934) all primitive Pythagorean triples can be generated using a set of linear transformations but it requires the triple (3,4,5) to start with. So we theorise a direct method to generate all possible primitive and non-primitive triples for a given number (one leg of a right angle triangle). Our strategy will be the observation that the difference between $c$ and $b$ (or between $c$ and $a$) can have only certain distinct values depending on the given number $a$ (or $b$).

Let us consider $c = b + \Delta$. So that equation (1) becomes
$$a^2 + b^2 = (b + \Delta)^2$$
$$\Rightarrow a^2 = \Delta(2b + \Delta) \qquad (2)$$
$$\Rightarrow b = \frac{a^2 - \Delta^2}{2\Delta} \qquad (3)$$

Equation (3) clearly shows that $\Delta$ must be a factor of $a^2$ for integral values of $b$. This is the first constraint that prevents $\Delta$ possessing any arbitrary value.

Further we note that $b$ to have positive value

$$a^2 - \Delta^2 > 0$$

$$\Rightarrow a > \Delta \tag{4}$$

This is the second condition imparted on $\Delta$.

In sections II and III we describe the case of Pythagorean triples extensively, in section IV we extend the theory for quadruples and finally in section V we discuss the generalization to n-tuples. Section VI provides a simple technique to generate tuples of arbitrary length starting from any single number.

## II. Generating Primitive Pythagorean Triple

A Pythagorean triple $(a, b, c)$ is said to be primitive when $a, b, c$ are coprime to each other i.e. $\gcd(a, b, c) = 1$. In this section we discuss for a given value of $a$ what are the possible values of $\Delta$. Before we proceed let us recall some basic properties of primitive triple. We know $c$ is always odd and only one of $a$ and $b$ is odd, so if we choose $a$ as even $\Delta$ will be even and if we set $a$ as odd $\Delta$ will also be odd.

Now numbers can be divided into three categories on the basis of their prime factorization.
(A) even numbers which are only powers of 2.
(B) odd numbers consisting of powers of any prime.
(C) even numbers containing both powers of 2 and other primes.

So we divide our discussion into three subsections:

(A) Let $a = 2^m$ and $\Delta = 2^r$, where $m$ and $r$ are whole numbers.

Equation (4) demands that $m > r$. Now from equation (2) we get

$$(2^m)^2 = 2^r(2b + 2^r)$$

$$\Rightarrow b = 2^{r-1}(2^{2m-2r} - 1) \tag{5}$$

Since in this case $a$ is even $b$ must be odd. Now the right hand side of equation (5) will be odd only if $r = 1$.
Hence the triple will be $(2^m, 2^{2m-2} - 1, 2^{2m-2} + 1)$.

*Example:* Let $a = 2^5 = 32$. So $\Delta = 2$.
Hence, $b = 255$ and $c = 257$. Thus we get (32,255,257) which is a primitive triple.
But if we choose $\Delta \neq 2$, say $\Delta = 8$ we obtain $b = 60$ and $c = 68$. Clearly the triple (32,60,68) is non-primitive that can be obtained from primitive triple (8,15,17).
So, in this case the only possible value of $\Delta$ is 2.

(B) Let $a = p^s k$, where $p$ is a prime factor of $a$ and $k$ is the product of other prime factors and $s$ can have any positive integral value. Then $\Delta = p^t$ where $t$ is also a whole number. From equation (2), we get

$$(p^s k)^2 = p^t(2b + p^t)$$

$$\Rightarrow b = \frac{1}{2}(p^{2s-t} k^2 - p^t) \tag{6}$$

Now for the triple to be primitive $b$ must not have $p$ as a factor. So, from this equation it is clear that either $t = 0$ or $t = 2s$ provided $a > \Delta$.

Hence the triple will be $(p^s k, \frac{k^2 - p^{2s}}{2}, \frac{k^2 + p^{2s}}{2})$ for $\Delta = p^{2s}$ or $(p^s k, \frac{p^{2s} k^2 - 1}{2}, \frac{p^{2s} k^2 + 1}{2})$ for $\Delta = 1$.

Here $\Delta = 1$ ensures us that for each odd number there exists at least one primitive triple of the form $(a, \frac{a^2 - 1}{2}, \frac{a^2 + 1}{2})$.

*Example:* Let us consider $a = 3^2 \times 11 = 99$. Then $\Delta$ can have two values, $\Delta = 3^4 = 81$ or $\Delta = 1$.
For $\Delta = 81$, $b = 20$, $c = 101$. Clearly (99,20,101) forms a primitive triple.
Next for $\Delta = 1$, $b = 4900$, $c = 4901$ which again makes a primitive triple (99,4900,4901).
In this case we also observe that if we take $p = 11$, we should have $\Delta = 11^2 = 121$ but it is not possible because it violates condition (4).

(C) Let $a = 2^m p^s k$, where $p$ is one prime factor of $a$ and $k$ is the product of other prime factors and $m, s$ can have any positive integral values. Then $\Delta = 2^r p^t$, where $r, t$ are whole numbers. From equation (2) we get

$$(2^m p^s k)^2 = 2^r p^t (2b + 2^r p^t)$$

$$\Rightarrow b = 2^{r-1}(2^{2m-2r} p^{2s-t} k^2 - p^t) \tag{7}$$

(i) When $m > r$, the terms within the parentheses on right hand side of equation (7) give odd number so that $b$ will be odd only if $r = 1$ and also $p$ must not be a factor of $b$ for the triple being primitive that requires $t = 0$ or $t = 2s$.

So, the primitive triple will be $(2^m p^s k, 2^{2m-2} p^{2s} k^2 - 1, 2^{2m-2} p^{2s} k^2 + 1)$ for $\Delta = 2$ and $(2^m p^s k, 2^{2m-2} k^2 - p^{2s}, 2^{2m-2} k^2 + p^{2s})$ for $\Delta = 2p^{2s}$.

*Example:* Let us consider $a = 2^3 5^2 \cdot 27 = 5400$. Then $\Delta$ can have two values, $\Delta = 2$ and $\Delta = 2 \cdot 5^4 = 1250$.
For $\Delta = 2$, $b = 7289999$, $c = 7290001$. So the primitive triple is (5400,7289999,7290001).
For $\Delta = 1250$, $b = 11039$, $c = 12289$. This gives (5400,11039,12289) which is also primitive.

(ii) When $m < r$, equation (7) becomes

$$b = 2^{r-1}(\frac{p^{2s-t} k^2}{2^{2r-2m}} - 1)$$

$$= 2^{2m-r-1}(p^{2s-t} k^2 - 2^{2r-2m} p^t) \tag{8}$$

Again the terms within the parentheses on right hand side of equation (8) give odd number so that $b$ will be odd only if $r = 2m - 1$ and also $b$ should not have $p$ as a factor for the triple being primitive which requires $t = 0$ or $t = 2s$.

*Example:* Let us consider $a = 2^3 3^2 \cdot 49 = 3528$. Then $\Delta$ will have two values $\Delta = 2^5 = 32$ and $\Delta = 2^5 3^4 = 2592$.
For $\Delta = 32$, $b = 194465$, $c = 194497$. So we get (3528,194465,194497) which is a primitive triple.
For $\Delta = 2592$, $b = 1105$, $c = 3697$. This gives (3528,1105,3697) which is also primitive.

In this case if we choose $p = 7$ instead of $p = 3$ it would give $\Delta = 2^5 7^4 = 76832$ which is not permissible since $\Delta$ becomes greater than $a$.

(iii) When $m = r$, from equation (7)

$$b = 2^{r-1}(p^{2s-t}k^2 - p^t) \qquad (9)$$

If $r \geq 1$, $b$ will always be even which gives non-primitive solutions. The only primitive triples are obtained when $r = m = 0$ but this will lead to equation (6) which we discuss earlier.

Here one important fact to be noted is that when $r = m = 1$, we never obtain a primitive solution so that integers of the configuration $(4n + 2)$ will always give non-primitive triple.

*Example:* The numbers 6, 10, 14, 18 ... will always form non-primitive triple.

$$6^2 + 8^2 = 10^2$$
$$10^2 + 24^2 = 26^2 \text{ etc.}$$

If we represent the given number as $a = 2^m p_1^{s_1} p_2^{s_2} \ldots p_n^{s_n}$ then $\Delta$ will be of the form $2^r q$ where $q = \prod_{i=1}^{N} q_i^{t_i}$ with $t_i = 0$ or $t_i = 2s_i$, $i = 1, 2, \ldots, n$.

*A general example:* Let us consider $a = 2^3 \times 3^2 \times 11 = 792$. Here $p_1 = 3$ and $p_2 = 11$. We show the different cases in the following table.

**Table 1:** Primitive triples for $a = 792$

| $\Delta$ | $b$ | $(a, b, c)$ |
| --- | --- | --- |
| $2^1 3^0 11^0 = 2$ | 156815 | (792,156815,156817) |
| $2^1 3^4 11^0 = 162$ | 1855 | (792,1855,2017) |
| $2^1 3^0 11^2 = 242$ | 1175 | (792,1175,1417) |
| $2^5 3^0 11^0 = 32$ | 9785 | (792,9785,9817) |
| $2^5 3^4 11^0 = 2592$ | ----- | ----- |
| $2^5 3^0 11^2 = 3872$ | ----- | ----- |
| $2^1 3^4 11^2 = 19602$ | ----- | ----- |
| $2^5 3^4 11^2 = 313632$ | ----- | ----- |

Here the last four values of $\Delta$ are not possible because of the condition (4). So only four primitive triples can be generated for $a = 792$.

Thus following the previous rules all possible primitive triples for a given number can be generated.

### III. Generating Non-primitive Pythagorean Triple

In the process of generating non-primitive Pythagorean triple the only constraint is $a > \Delta$. So we first need to factorize the given number and then $\Delta$ will be any combination of those factors except the cases for primitive triples. So it is obvious that if $a$ is even $\Delta$ must be even and if $a$ is odd $\Delta$ will also be odd. We illustrate the method by the following example.

Let $a = 60 = 2^2 \times 3 \times 5$. Again we show various cases in the following table.

**Table 2:** Non-primitive triples for $a = 60$

| $\Delta$ | $b$ | $(a, b, c)$ |
|---|---|---|
| $2^2 3^0 5^0 = 4$ | 448 | (60,448,452) |
| $2^1 3^1 5^0 = 6$ | 297 | (60,297,303) |
| $2^1 3^0 5^1 = 10$ | 175 | (60,175,185) |
| $2^2 3^1 5^0 = 12$ | 144 | (60,144,156) |
| $2^2 3^0 5^1 = 20$ | 80 | (60,80,100) |
| $2^3 3^1 5^0 = 24$ | 63 | (60,63,87) |
| $2^1 3^1 5^1 = 30$ | 45 | (60,45,75) |
| $2^2 3^2 5^0 = 36$ | 32 | (60,32,68) |
| $2^3 3^0 5^1 = 40$ | 25 | (60,25,65) |

Here $\Delta = 2, 8, 18$ and $50$ have not been taken because they generate primitive triples and the other combinations of the factors have been discarded as they violate condition (4).

Thus finding all the possible values of $\Delta$ we can obtain all non-primitive triples.

## IV. Pythagorean Quadruple

A Pythagorean quadruple is an ordered quadruple of positive integers $(a, b, c, d)$ such that

$$a^2 + b^2 + c^2 = d^2 \tag{10}$$

In this section we discuss how to generate all possible Pythagorean quadruples for a given set of $(a, b)$. Let $a^2 + b^2 = k$ and $d = c + \Delta$. Then equation (10) becomes

$$k + c^2 = (c + \Delta)^2$$

$$\Rightarrow c = \frac{k - \Delta^2}{2\Delta} \tag{11}$$

Three facts are clear from equation (11):

(i)  if $k$ is even $\Delta$ must be even and if $k$ is odd $\Delta$ must be odd for integral value of $c$.
(ii) when $k$ is even it ought to be an integral multiple of $2\Delta$.
(iii) for positive value of $c$ it is required that
$$k > \Delta^2 \tag{12}$$

We are mainly interested in the generation of primitive quadruples which we will discuss in three sections.

**Case (A): $a$ is even and $b$ is odd (or $b$ is even and $a$ is odd)**

Here $k$ is an odd number, so $\Delta$ is odd.

(i) We first consider that $a$ and $b$ have common factors $p_1, p_2, \ldots p_n$, then $k$ can be represented as $k = p_1^{m_1} p_2^{m_2} \ldots p_n^{m_n} q_1^{s_1} q_2^{s_2} \ldots q_N^{s_N}$ and $\Delta$ will be of the form $\Delta = p_1^{r_1} p_2^{r_2} \ldots p_n^{r_n} q_1^{t_1} q_2^{t_2} \ldots q_N^{t_N}$; $m_i, s_i, r_i$ and $t_i$ all being integers for all $i$. Equation (11) then gives

$$c = \frac{1}{2}(p_1^{m_1-r_1} p_2^{m_2-r_2} \ldots p_n^{m_n-r_n} q_1^{s_1-t_1} q_2^{s_2-t_2} \ldots q_N^{s_N-t_N} - p_1^{r_1} p_2^{r_2} \ldots p_n^{r_n} q_1^{t_1} q_2^{t_2} \ldots q_N^{t_N}) \tag{13}$$

So for primitive solution $r_i = 0$ or $r_i = m_i$, $i = 1, 2, \ldots, n$; and $t_j$ can take all integral values from 0 to $s_j$ with the restriction given by equation (12). Then $\Delta = \prod_{i=1}^{n} p_i^{r_i} \prod_{j=1}^{N} q_j^{t_j}$, with either $r_i = 0$ or $r_i = m_i$ for all $i$ and $t_j$ has those values discussed above.

*Example:* (1) $a = 12$ and $b = 15$. Then $k = 369 = 3^2 \times 41$.
So $\Delta = 1, 3^2$ but $\Delta = 41$ is not possible as $\Delta^2 > k$.
Thus we get,
$$c = 184 \text{ and } d = 185 \text{ when } \Delta = 1$$
$$c = 16 \text{ and } d = 25 \text{ when } \Delta = 3^2$$
So the primitive quadruples for $a = 12$ and $b = 15$ are (12,15,184,185) and (12,15,16,25).

(2) Let $a = 2 \times 3 \times 5 \times 7 = 210$ and $b = 3^3 \times 5 = 135$. Then $k = 3^2 \times 5^2 \times 277 = 62325$
So $\Delta = 1, 3^2, 5^2, 3^2 \times 5^2$. Other combinations of $\Delta$ are not possible due to equation (12).

**Table 3:** Primitive Pythagorean quadruples for $a = 210$ and $b = 135$

| Δ | c | d | $(a,b,c,d)$ |
|---|---|---|---|
| $3^0 5^0 277^0 = 1$ | 31162 | 31163 | (210,135,31162,31163) |
| $3^2 5^0 277^0 = 9$ | 3458 | 3467 | (210,135,3458,3467) |
| $3^0 5^2 277^0 = 25$ | 1234 | 1259 | (210,135,1234,1259) |
| $3^2 5^2 277^0 = 225$ | 26 | 251 | (210,135,26,251) |

(ii) Now we consider that $a$ and $b$ have no common factors so that $k = q_1^{s_1} q_2^{s_2} \ldots q_N^{s_N}$, then $\Delta = \prod_{j=1}^{N} q_j^{t_j}$ where $t_j$ can take all integral values from 0 to $s_j$ provided $k > \Delta^2$.

*Example:* Let $a = 8$ and $b = 19$. Then $k = 425 = 5^2 \times 17$. ∴ $\Delta = 1, 5, 17$.

**Table 4:** Primitive Pythagorean quadruples for $a = 8$ and $b = 19$

| Δ | c | d | $(a,b,c,d)$ |
|---|---|---|---|
| $5^0 17^0 = 1$ | 212 | 213 | (8,19,212,213) |
| $5^1 17^0 = 5$ | 40 | 45 | (8,19,40,45) |
| $5^0 17^1 = 17$ | 4 | 21 | (8,9,4,21) |

Here we draw an important conclusion that whenever one of $a$ or $b$ is odd, there will be at least one primitive quadruple with $\Delta = 1$.

### Case (B): Both $a$ and $b$ are even

Here $k$ is an even number, so $\Delta$ will be even. This case will be same as discussed in case (A) except the introduction of some power of 2 so that now $k = 2^m p_1^{m_1} p_2^{m_2} \ldots p_n^{m_n} q_1^{s_1} q_2^{s_2} \ldots q_N^{s_N}$ and $\Delta = 2^r p_1^{r_1} p_2^{r_2} \ldots p_n^{r_n} q_1^{t_1} q_2^{t_2} \ldots q_N^{t_N}$. Equation (11) then gives

$$c = 2^{m-r-1} p_1^{m_1-r_1} p_2^{m_2-r_2} \ldots p_n^{m_n-r_n} q_1^{s_1-t_1} q_2^{s_2-t_2} \ldots q_N^{s_N-t_N} - 2^{r-1} p_1^{r_1} p_2^{r_2} \ldots p_n^{r_n} q_1^{t_1} q_2^{t_2} \ldots q_N^{t_N} \quad (14)$$

So, the conditions for obtaining primitive quadruples are same with an additional condition either $r = 1$ or $r = m - 1$. Then the values of $\Delta$ are

$\Delta = 2^1 \prod_{i=1}^{n} p_i^{r_i} \prod_{j=1}^{N} q_j^{t_j}$ or $\Delta = 2^{m-1} \prod_{i=1}^{n} p_i^{r_i} \prod_{j=1}^{N} q_j^{t_j}$, here $r_i = 0$ or $n_i$ and $t_j$ has any integral value from 0 to $s_j$, provided $k > \Delta^2$.

Again an important conclusion from this theory is that for any two given even numbers there always exists a primitive quadruple with $\Delta = 2$.

*Example:* Let $a = 2 \times 3 = 6$ and $b = 2 \times 3 \times 5 = 30$. Then $k = 2^3 \times 3^2 \times 13 = 936$. The following table shows all possible primitive cases.

**Table 5:** Primitive Pythagorean quadruples for $a = 6$ and $b = 30$

| Δ | c | d | (a, b, c, d) |
|---|---|---|---|
| $2^1 3^0 13^0 = 2$ | 233 | 235 | (6,30,233,235) |
| $2^2 3^0 13^0 = 4$ | 115 | 119 | (6,30,115,119) |
| $2^1 3^2 13^0 = 18$ | 17 | 35 | (6,30,17,35) |
| $2^1 3^0 13^1 = 26$ | 5 | 31 | (6,30,5,31) |

$\Delta \neq 2^2 \times 3^2, 2 \times 3^2 \times 13, 2^2 \times 13, 2^2 \times 3^2 \times 13$, as their squares are greater than $k$.

### Case (C): Both $a$ and $b$ are odd

To study this case we represent $a$ and $b$ as $a = 2x + 1$ and $b = 2y + 1$ where $x, y$ are positive integers. Then

$$k = a^2 + b^2 = (2x + 1)^2 + (2y + 1)^2$$
$$\Rightarrow k = 4(x^2 + y^2) + 4(x + y) + 2 \tag{15}$$

Since $k$ is even $\Delta$ must be even but then $k$ must contain a factor 4 for being an integral multiple of $2\Delta$. But equation (15) clearly shows that $k$ is a multiple of 2 only i.e. $k = 2 \times odd\ number$.

So no set of quadruples (primitive and non-primitive) can be obtained if both the given numbers are odd.

To produce non-primitive quadruples, we first calculate $k$ and factorise it. Then $\Delta$ can have any combination of those factors expect those for primitive cases with the condition $k > \Delta^2$. Following examples will make it clear.

*Example*: (i) Let $a = 105 = 3 \times 5 \times 7$ and $b = 150 = 3 \times 5^2 \times 2$.
Then $k = 3^2 \times 5^2 \times 149 = 33525$. So the values of $\Delta$ that give non-primitive quadruples will be $\Delta = 3, 5, 3 \times 5, 3 \times 5^2$.

**Table 6:** Non-primitive Pythagorean quadruples for $a = 105$ and $b = 150$

| Δ | c | d | (a, b, c, d) |
|---|---|---|---|
| $3^1 5^0 149^0 = 3$ | 5586 | 5589 | (105,150,5586,5589) |
| $3^0 5^1 149^0 = 5$ | 3350 | 3355 | (105,150,3350,3355) |
| $3^1 5^1 149^0 = 15$ | 1110 | 1125 | (105,150,1110,1125) |
| $3^2 5^1 149^0 = 45$ | 350 | 395 | (105,150,350,395) |
| $3^1 5^2 149^0 = 75$ | 186 | 261 | (105,150,186,261) |

Here $\Delta = 1, 149, 3^2, 5^2$ give primitive quadruples and the rest combinations violate condition (12), so they are not taken.

(ii) Let $a = 2 \times 7 = 14$ and $a = 2 \times 7^2 = 98$. The $k = 2^3 \times 7^2 \times 5^2 = 9800$. Here the primitive cases are

**Table 7:** Non-primitive Pythagorean quadruples for $a = 14$ and $b = 98$

| Δ | c | d | (a, b, c, d) |
|---|---|---|---|
| $2^1 7^0 5^1 = 10$ | 485 | 495 | (14,98,485,495) |
| $2^1 7^1 5^0 = 14$ | 343 | 357 | (14,98,343,357) |
| $2^2 7^0 5^1 = 20$ | 235 | 255 | (14,98,235,255) |
| $2^2 7^1 5^0 = 28$ | 161 | 189 | (14,98,161,189) |
| $2^1 7^1 5^1 = 70$ | 35 | 105 | (14,98,35,105) |

Here Δ= $2, 2 \times 7^2, 2 \times 5^2, 2^2$ give primitive sets and the rest do not obey the condition (12). So these are not taken.

In this way all Pythagorean Quadruples can be generated for any two given numbers.

## V. Pythagorean n-tuple

A Pythagorean n-tuple is a set of n positive integers $(a_1, a_2, \ldots, a_n)$ such that

$$a_1^2 + a_2^2 + \cdots + a_{n-1}^2 = a_n^2 \tag{16}$$

When $(n - 2)$ numbers are given, we can calculate the rest two and form an n-tuple $(a_1, a_2, \ldots, a_n)$. This process is quite similar to the process or generating primitive quadruples with $a_1^2 + a_2^2 + \cdots + a_{n-2}^2 = k$ and $a_n - a_{n-1} = \Delta$.

Now among given $(n - 2)$ numbers some will be odd, say $\lambda$ numbers and the rest will be even.

(i) When $\lambda = 2$, $k = 2 \times odd\ number$. This does not lead to the formation of n-tuple which is similar to case (C) of quadruple.

(ii) When $\lambda = odd$, $k$ is odd too and $k = p_1^{m_1} p_2^{m_2} \ldots p_n^{m_n} q_1^{s_1} q_2^{s_2} \ldots q_N^{s_N}$, where $p_1, p_2, \ldots, p_n$ with some powers are common factors of the given $(n - 2)$ integers. Then Δ= $p_1^{r_1} p_2^{r_2} \ldots p_n^{r_n} q_1^{t_1} q_2^{t_2} \ldots q_N^{t_N}$, $r_i = 0$ or $n_i$ and $t_j$ can be any integer from 0 to $s_j$ satisfying $k > \Delta^2$.

*Example:* Let $a_1 = 55, a_2 = 15, a_3 = 20, a_4 = 10, a_5 = 35, a_6 = 45, a_7 = 30$ and $a_8 = 25$.
∴ $k = 5^2 \times 11 \times 31 = 8525$ and Δ= $1, 5^2, 11, 31$.

**Table 8:** Primitive Pythagorean n-tuples (n=10) for 8 given numbers

| Δ | $a_9$ | $a_{10}$ |
|---|---|---|
| $5^0 11^0 31^0 = 1$ | 4262 | 4263 |
| $5^0 11^1 31^0 = 11$ | 382 | 393 |
| $5^2 11^0 31^0 = 25$ | 158 | 183 |
| $5^0 11^0 31^1 = 31$ | 122 | 153 |

So the primitive n-tuples (n=10) for those 8 numbers are

$55^2 + 15^2 + 20^2 + 10^2 + 35^2 + 45^2 + 30^2 + 25^2 + 4262^2 = 4263^2$

$55^2 + 15^2 + 20^2 + 10^2 + 35^2 + 45^2 + 30^2 + 25^2 + 382^2 = 393^2$

$55^2 + 15^2 + 20^2 + 10^2 + 35^2 + 45^2 + 30^2 + 25^2 + 158^2 = 183^2$

$55^2 + 15^2 + 20^2 + 10^2 + 35^2 + 45^2 + 30^2 + 25^2 + 122^2 = 153^2$

(iii) When $m_1 = even$ (greater than 2) or $m_1 = 0$ i.e. all the given numbers are even, then $k$ is also even and it will have the form $k = 2^m p_1^{n_2} p_2^{n_2} \dots p_n^{n_n} q_1^{s_1} q_2^{s_2} \dots q_1^{s_N}$, $m > 1$, where $p_1, p_2, \dots, p_n$ or some powers of them are common factors of the given numbers and $q_1, q_2, \dots, q_N$ are other primes. Then $\Delta = 2^r \prod_{i=1}^{n} p_i^{r_i} \prod_{j=1}^{N} q_j^{t_j}$, $r = 1$ or $(m-1)$, $r_i = 0$ or $n_i$ and $t_j$ can take any integral value from 0 to $s_j$ provided $k > \Delta^2$.

*Example:* Let $a_1 = 24$, $a_2 = 57$, $a_3 = 54$, $a_4 = 33$, $a_5 = 39$, $a_6 = 21$ and $a_7 = 48$. Now $k = 3^3 \times 2^6 \times 7 = 12096$. So $\Delta = 2^1 \times 3^0, 2^1 \times 3^3, 2^1 \times 3^0 \times 7, 2^5 \times 3^0$ for primitive solutions.

**Table 9:** Primitive Pythagorean n-tuples (n=9) for 7 given numbers

| Δ | $a_8$ | $a_9$ |
|---|---|---|
| $2^1 3^0 7^0 = 2$ | 3023 | 3025 |
| $2^1 3^0 7^0 = 14$ | 425 | 439 |
| $2^5 3^0 7^0 = 32$ | 173 | 205 |
| $2^1 3^3 7^0 = 54$ | 85 | 139 |

So the primitive n-tuples ($n = 9$) for the given 7 numbers.

$24^2 + 57^2 + 54^2 + 33^2 + 39^2 + 21^2 + 48^2 + 3023^2 = 3025^2$

$24^2 + 57^2 + 54^2 + 33^2 + 39^2 + 21^2 + 48^2 + 425^2 = 439^2$

$24^2 + 57^2 + 54^2 + 33^2 + 39^2 + 21^2 + 48^2 + 173^2 = 205^2$

$$24^2 + 57^2 + 54^2 + 33^2 + 39^2 + 21^2 + 48^2 + 85^2 = 139^2$$

We can generate the non-primitive n-tuples also in a manner similar to the case of quadruple discussed previously.

## VI. Generating Pythagorean n-tuple starting from a single number

Having the discussion of generating different Pythagorean tuples elaborately here we give a simple method to generate tuples of arbitrary length starting from a given number using the theory of Pythagorean triple.

If the number $a_1$ is given, we can calculate the triple $(a_1, a_2, b_3)$ i.e. $a_1^2 + a_2^2 = b_3^2$. Then starting from $b_3$, $(b_3, a_3, b_4)$ can be obtained. So the equation extends to $a_1^2 + a_2^2 + a_3^2 = b_4^2$, a quadruple. Proceeding in the same way we can elongate the chain and after $(n-2)$ iterations we obtain the n-tuple:
$$a_1^2 + a_2^2 + a_3^2 + \cdots + a_{n-1}^2 = b_n^2 = a_n^2.$$

Now it is to be noted that for a given number $(a_1)$ we can have several triples, so each of them will form a different branch and each branch will have a sub-branch and so on. A branch will be primitive only if gcd of any three numbers is 1. In this way we can obtain several n-tuples from a given number but this method is incapable of generating all possible cases since the method totally depends on the process of generating triples.

*Example:* Here we show three branches with $a = 15$.
$15^2 + 36^2 + 760^2 + 289560^2 = 289561^2$
$15^2 + 20^2 + 60^2 + 2112^2 + 2232384^2 = 2232385^2$
$15^2 + 8^2 + 144^2 + 348^2 + 71064^2 = 71065^2$

## VII. Conclusion

The major advantage of our method is that it does not require any primitive set to start with and finding proper multipliers or transformations to obtain the desired tuple. An interesting fact is that just by factorizing $k$ we can forecast how many primitive and non-primitive cases are possible before actually calculating them which cannot be done by classical formulas. We are also able to produce tuples of any length starting from a given number in a very simple way.